\documentclass{amsart} 

\usepackage{amsmath,amssymb,latexsym, amscd}
\usepackage{exscale, cite, eps fig, graphics}

\renewcommand{\geq}{\geqslant}
\renewcommand{\leq}{\leqslant}

\newtheorem*{corollary*}{Corollary}

\newtheorem*{theorem*}{Theorem}
\newtheorem*{remark*}{Remark}
\numberwithin{equation}{section}

\title[Solitary waves in deep water]
{Kinetic, potential and surface tension energies \\ of solitary waves in deep water}

\author[Hur]{Vera~Mikyoung~Hur}
\address{Department of Mathematics, University of Illinois at Urbana-Champaign, Urbana, IL 61801 USA}
\email{verahur@math.uiuc.edu}

\date{\today}


\begin{document}

\maketitle

\begin{abstract}
We present an exact relation among the kinetic, potential and surface tension energies 
of a solitary wave in deep water in all dimensions. 
We deduce its non-existence in the absence of the effects of surface tension, 
provided that gravity acts in a direction opposite to what is physically realistic. 
\end{abstract}

\section{Introduction}\label{sec:intro}

It is a matter of experience that waves which are commonly seen in the ocean or a lake propagate 
in a certain direction approximately at a constant speed without change of form, namely traveling waves. 
They have historically stimulated a considerable part of the development in the theory of wave motion,
from Stokes' conjecture about the wave of maximum amplitude
to Russell's famous horseback observations and to the elucidation of the Korteweg-de Vries solitons.

By a {\em solitary water wave}, we mean a traveling wave solution, 
for which the fluid surface asymptotically approaches a constant level over a nearly shear flow. 
In two dimensions in the finite depth case, 
their rigorous mathematical theory dates back to 
the constructions in \cite{FH54} and \cite{Beale77} of small amplitude waves,
and it includes the global bifurcation result in \cite{AT-solitary}, 
extensions in \cite{AKcapillary}, \cite{Beale-capillary}, \cite{Sun-capillary} and
\cite{IKsumm}, \cite{IK}, \cite{BGT},
among others, to capillary-gravity waves, 
and in \cite{Hur-solitary}, \cite{GW08} and \cite{Wheeler1} to waves with vorticity. 
The symmetry and monotonicity properties were discussed in \cite{CSsymm} and \cite{Hur-symmS},
and the regularity properties were in \cite{Lewy} and \cite{Hur-anal}. 
Non-uniqueness and linear instability were addressed in \cite{Plotnikov} and in \cite{Lin-ww}. 

\

In stark contrast, Sun in \cite{Sun97} argued for the non-existence of solitary waves 
to interfacial fluids problems in two dimensions in the infinite depth case, provided that 
the profile is either positively elevated above the mean fluid level or negatively depressed.
Craig in \cite{Craig02} made another proof 
for gravity driven and positively elevated, solitary water waves 
in two dimensions in the infinite depth case and in three dimension in the finite depth case.
Therefore one concludes that a solitary water wave 
in two dimensions in the infinite depth case and in three dimensions in the finite depth case
cannot be everywhere positive (or negative, either, in the former setting). 
As Craig emphasized in \cite{Craig02}, 
it ``leaves open the tantalizing possibility" of solitary water waves, which oscillate about the mean fluid level.
In the presence of the effects of surface tension, solitary water waves 
indeed arise in two dimensions in the infinite depth case (see \cite{IKirrmann}, for instance)
and in three dimensions in the finite depth case (see \cite{GS08}, for instance), 
which necessarily change sign. In three dimensions in the finite depth case, 
they are found near the Kadomtsev-Petviashvili (KP) solitary waves. 
In the absence of the effects of surface tension, however,
the KP equation does not invite a localized steady solution, 
bearing out that gravity driven, solitary water waves may not exist in three dimensions.
Moreover small-data global-existence results in \cite{Wu3g} and \cite{GMSg} for the initial value problem
indicate that gravity-driven solitary water waves of small amplitudes are unlikely to exist 
in three dimensions in the infinite depth case.

Recently in \cite{Hur-no1}, the author eliminated the positivity requirement 
from \cite{Craig02} and \cite{Sun97} in two dimensions in the infinite depth case. 
The proof relates the free boundary problem in potential theory with a nonlinear pseudodifferential equation
via conformal mapping techniques and derives a Pohozaev type identity.
Conformal mappings are not available in higher dimensions.
Nevertheless, here we apply the Pohozaev identity technique 
to the free boundary problem (see \eqref{E:potential})
and we find an exact relation among the kinetic, potential and surface tension energies of 
a solitary wave in deep water in all dimensions. 
In two dimensions we provide another proof of an integral identity in \cite{Sun97}.

We conclude that in the absence of the effects of surface tension,
the solitary wave problem in deep water admits no nontrivial solutions,
provided that gravity acts in a direction opposite to what is physically realistic. 
The result partly recovers that in \cite{Hur-no1}.

\

Other exact relations among integral quantities were known 
for solitary waves in two dimensions in the finite depth case; 
see \cite{LH-relations}, \cite{McCowan}, \cite{Starr}, among others. 
Similar relations would perhaps be useful to improve the conclusion herein, 
but the proofs seem to break down in the infinite depth case.

\section{Results}\label{sec:result}

The water wave problem in the simplest form concerns 
the wave motion at the interface separating an incompressible inviscid fluid below a body of air,
acted upon by gravity and possibly surface tension.
The flow in the bulk of the fluid satisfies the Euler equations in hydrodynamics with the force of gravity.
The kinematic and dynamic conditions at the fluid surface state, respectively, that 
fluid particles do not invade the air, nor vice versa, 
and that the jump in pressure across the surface, in the presence of the effects of surface tension,
is proportional to the mean curvature. 
The flow is assumed to be nearly at rest at great depths. 
For definiteness, we shall work in the $(\mathbf{x},y)$-coordinates, 
where $\mathbf{x}\in\mathbb{R}^{n-1}$, $n\geq 2$ an integer, is horizontal and $y \in \mathbb{R}$ is vertical. 
If the fluid surface is in the graph form $y=\eta(\mathbf{x},t)$, say, 
and if the flow beneath it is irrotational then the governing equations may be expressed 
in terms of the velocity potential $\phi(\mathbf{x},y,t)$ as
\begin{equation}\label{E:WW}
\begin{aligned}
&\Delta\phi=0 &\text{in}\quad&-\infty<y<\eta(\mathbf{x},t),\\
&\eta_t+\nabla_\mathbf{x}\phi\cdot\nabla\eta-\phi_y=0 &\text{at}\quad&y=\eta(\mathbf{x},t),\\
&\phi_t+\frac12|\nabla\phi|^2+g\eta-T\nabla\Big(\frac{\nabla\eta}{\sqrt{1+|\nabla\eta|^2}}\Big)=0
&\text{at}\quad&y=\eta(\mathbf{x},t),\\
&\nabla\phi\to\mathbf{0}&\text{as}\quad &y\to-\infty.
\end{aligned}
\end{equation}
Here $g$ denotes the constant due to gravitational acceleration
and $T\geq 0$ is the coefficient of surface tension. Note that $g>0$ is physically realistic.
In the case of $g<0$, the initial value problem associated with the linear part of \eqref{E:WW} is ill-posed;
indeed plane waves of the form $e^{i(\mathbf{k}\cdot\mathbf{x}-\omega t)}$ obey 
the dispersion relation $\omega^2=g|\mathbf{k}|$, 
which permits unbounded growths for high frequency waves. 
But {\em existence}, and not ``stability", is the issue here.
Throughout subscripts denote partial differentiation;
$\nabla$ means the usual gradient in $\mathbb{R}^n$ or $\mathbb{R}^{n-1}$,
and $\nabla_\mathbf{x}$ is the gradient in $\mathbb{R}^{n-1}$ in the horizontal variables.

\

We may define the kinetic, potential and surface tension energies of a solution to \eqref{E:WW} as
\[
\int^{\eta(\mathbf{x},t)}_{-\infty}\int_{\mathbb{R}^{n-1}} \frac12|\nabla\phi|^2~d\mathbf{x}dy,\qquad
\int_{\mathbb{R}^{n-1}}\frac12g\eta^2~d\mathbf{x}\quad\text{and}\quad
\int_{\mathbb{R}^{n-1}}\frac12T(\sqrt{1+|\nabla\eta|^2}-1)~d\mathbf{x},
\]
respectively. Their sum is the total energy and it is a constant of motion; 
see \cite{Zakharov-WW}, for instance.

\

The solitary wave problem seeks for a solution, for which 
the fluid surface and the velocity potential depend upon $(\mathbf{x}-\mathbf{c}t, y)$, 
where $\mathbf{c}\in\mathbb{R}^{n-1}$ is the speed of wave propagation,
the fluid surface asymptotically tends to zero, say, and the flow in the far field is nearly uniform, i.e. 
\[
\eta(\mathbf{x})\to0 \quad\text{as $|\mathbf{x}|\to\infty$}\quad\text{and}\quad
\nabla\phi(\mathbf{x},y)\to \mathbf{0}\quad \text{as $|(\mathbf{x},y)|\to \infty$.}
\]
In the frame of reference moving at the velocity of wave propagation, let 
\[
\Omega=\{(\mathbf{x},y): \mathbf{x}\in\mathbb{R}^{n-1}, -\infty<y<\eta(\mathbf{x})\}
\quad \mbox{and}\quad
\Gamma=\{(\mathbf{x}, \eta(\mathbf{x})):  \mathbf{x}\in\mathbb{R}^{n-1}\}
\]
represent, respectively, the (stationary) fluid domain and its surface. 
If we furthermore assume that $\phi(\mathbf{x},y)\to0$ as $|(\mathbf{x},y)|\to \infty$ then 
the kinetic energy, after an application of the divergence theorem, becomes 
\[
\frac12\int_\Gamma \phi(\nabla\phi\cdot\boldsymbol{\nu})~dS,
\]
where 
\[
\boldsymbol{\nu}=\frac{1}{\sqrt{1+|\nabla\eta(\mathbf{x})|^2}}(-\nabla\eta,1)
\]
denotes the outward pointing unit normal vector along the fluid surface
and $dS$ is the surface measure of $\Gamma$. 
This allows us to handle various integrals in the following section.
In summary, the solitary wave problem in deep water reads: 

find $\eta$ defined in $\mathbb{R}^{n-1}$, 
$\phi$ defined over $\overline{\Omega}$ and a parameter $\mathbf{c}\in\mathbb{R}^{n-1}$ such that 
\begin{subequations}\label{E:potential}
\begin{alignat}{2}
&\Delta \phi=0 \qquad &&\mbox{in $\Omega$}, \label{E:laplace} \\
&\phi_y=(\nabla_\mathbf{x}\phi-\mathbf{c})\cdot\nabla\eta \qquad&&\mbox{at $\Gamma$},
\label{E:kinematic} \\
&\frac12|\nabla\phi|^2-\mathbf{c}\cdot\nabla\phi+gy=T\nabla\Big(\frac{\nabla\eta}{\sqrt{1+|\nabla\eta|^2}}\Big)
\qquad&&\mbox{at $\Gamma$}, 
\label{E:bernoulli} \\
&\eta(\mathbf{x})\to 0 \qquad & &\mbox{as $|\mathbf{x}|\to \infty$,} \\
&\phi(\mathbf{x},y)\to 0\qquad&&\mbox{as $|(\mathbf{x},y)|\to \infty$.} \label{E:infinity}
\end{alignat}
\end{subequations}
In two dimensions, \eqref{E:potential} may be reformulated into a single, nonlinear pseudodifferential equation
via conformal mapping techniques; see \cite{Hur-no1}, for instance.

\

In the absence of the effects of surface tension, i.e. $T=0$, 
a straightforward calculation reveals that \eqref{E:potential} remains invariant under 
\[
\eta(\mathbf{x})\mapsto \lambda^{-1}\eta(\lambda\mathbf{x}), \qquad
\phi(\mathbf{x},y)\mapsto \lambda^{-3/2}\phi(\lambda\mathbf{x},\lambda y)
\quad\text{and}\quad \mathbf{c}\mapsto\lambda^{-1/2}\mathbf{c}
\]
for all $\lambda>0$. The vector field $(\mathbf{x}, y)\cdot\nabla$ generates the scaling symmetry,
and it suggests that we apply the Pohozaev identity technique to \eqref{E:potential}.
Incidentally the Pohozaev identity technique was devised in \cite{Pohozaev}
and it led to a wide range of applications; see \cite{PucciSerrin}, for instance.

\

In what follows,
\[
[f(\mathbf{x},y)]_\Gamma=f(\mathbf{x},\eta(\mathbf{x}))
\]
denotes $f$ restricted to $\Gamma$. We state the main result.

\begin{theorem*}
In the case of $T=0$, 
if $\eta\in C^1(\mathbb{R}^{n-1})$ and $\phi\in C^\infty(\Omega)\cap C(\overline{\Omega})$ 
solve \eqref{E:potential} and, moreover, if 
$\eta\in H^1(\mathbb{R}^{n-1})$, $\mathbf{x}\cdot\nabla\eta \in L^2(\mathbb{R}^{n-1})$ and
$\nabla\phi\in L^2(\Omega)$, $(\mathbf{x},y)\cdot\nabla\phi\in L^2(\Omega)$,
$[\phi]_\Gamma \in H^1(\mathbb{R}^{n-1})$, $\mathbf{x}\cdot[\nabla\phi]_\Gamma \in L^2(\mathbb{R}^{n-1})$,
then
\begin{equation}\label{E:pohozaev}
\frac{n}{2}\iint_\Omega |\nabla \phi|^2~d\mathbf{x}dy
=\frac{n+1}{2}g\int_{\mathbb{R}^{n-1}} \eta^2~d\mathbf{x}.
\end{equation}

In the case of $T\neq 0$, if in addition $\mathbf{x}\cdot \nabla^2\eta \in L^2(\mathbb{R}^{n-1})$ then
\begin{equation}\label{E:pohozaevT}
\frac{n}{2}\iint_\Omega |\nabla \phi|^2~d\mathbf{x}dy
=\frac{n+1}{2}g\int_{\mathbb{R}^{n-1}} \eta^2~d\mathbf{x}
+\frac12T\int_{\mathbb{R}^{n-1}}(\sqrt{1+|\nabla\eta|^2}-1)~d\mathbf{x}.
\end{equation}
\end{theorem*}

In two dimensions \eqref{E:pohozaevT} agrees with (1.3) in \cite{Sun97}.

\

\begin{corollary*}
In the case of $T=0$, 
assume that $\eta\in C^1(\mathbb{R}^{n-1})$ and $\phi\in C^\infty(\Omega)\cap C(\overline{\Omega})$ 
solve \eqref{E:potential} and, 
$\eta\in H^1(\mathbb{R}^{n-1})$, $\mathbf{x}\cdot\nabla\eta \in L^2(\mathbb{R}^{n-1})$ and
$\nabla\phi\in L^2(\Omega)$, $(\mathbf{x},y)\cdot\nabla\phi\in L^2(\Omega)$,
$[\phi]_\Gamma \in H^1(\mathbb{R}^{n-1})$, $\mathbf{x}\cdot[\nabla\phi]_\Gamma \in L^2(\mathbb{R}^{n-1})$.
Then $\eta\equiv0$ and $\phi\equiv0$ if $g\leq 0$.

\end{corollary*}

In two dimensions it partly recovers that in \cite{Hur-no1}.

\

\section{Proof of Theorem}\label{sec:proof}

We may assume that $\eta$ and $\phi$, up to boundary, are smooth and 
decay to zero at infinity faster than polynomials together with their derivatives of all orders.
We may mollify $\eta$ and $\phi$ and appeal to the dominated convergence theorem
to ensure that all the integrals below converge under the hypothesis of Theorem. 

\

Suppose for now that $T=0$.
We multiply \eqref{E:laplace} by $\phi$ and integrate over $\Omega$ to show that
\begin{align*}
0=\iint_\Omega \phi\Delta\phi~d\mathbf{x}dy
=&\int_\Gamma\phi(\nabla\phi\cdot\boldsymbol{\nu})~dS-\iint_\Omega|\nabla\phi|^2~d\mathbf{x}dy \\
=&\int_{\mathbb{R}^{n-1}}[\phi]_\Gamma
([\phi_y]_\Gamma-[\nabla_\mathbf{x}\phi]_\Gamma\cdot\nabla\eta)~d\mathbf{x}
-\iint_\Omega|\nabla\phi|^2~d\mathbf{x}dy  \\
=&-\int_{\mathbb{R}^{n-1}}[\phi]_\Gamma(\mathbf{c}\cdot\nabla)\eta~d\mathbf{x}
-\iint_\Omega|\nabla\phi|^2~d\mathbf{x}dy.
\end{align*}
Recall that 
\[
\boldsymbol{\nu}=\frac{1}{\sqrt{1+|\nabla\eta(\mathbf{x})|^2}}(-\nabla\eta,1)
\]
denotes the outward pointing unit normal vector along the fluid surface 
and $dS$ is the surface measure of $\Gamma$.
The first equality uses \eqref{E:laplace}, the second equality uses the divergence theorem, 
and the last equality uses \eqref{E:kinematic}. 
Indeed $[\phi_y]=([\nabla_\mathbf{x}\phi]-\mathbf{c})\cdot\nabla\eta$. 
Moreover an integration by parts leads to that
\begin{align*}
\int_{\mathbb{R}^{n-1}}[\phi]_\Gamma(\mathbf{c}\cdot\nabla)\eta~d\mathbf{x}
=-\int_{\mathbb{R}^{n-1}}\eta(\mathbf{c}\cdot\nabla)[\phi]_\Gamma~d\mathbf{x} 
=-\int_{\mathbb{R}^{n-1}}\eta\mathbf{c}\cdot
([\nabla_\mathbf{x}\phi]_\Gamma+[\phi_y]_\Gamma\nabla\eta)~d\mathbf{x}.
\end{align*}
Therefore
\begin{equation}\label{E:pohozaev1}
\iint_\Omega|\nabla\phi|^2~d\mathbf{x}dy
=-\int_{\mathbb{R}^{n-1}}[\phi]_\Gamma(\mathbf{c}\cdot\nabla)\eta~d\mathbf{x}
=\int_{\mathbb{R}^{n-1}}\eta\mathbf{c}\cdot
([\nabla_\mathbf{x}\phi]_\Gamma+[\phi_y]_\Gamma\nabla\eta)~d\mathbf{x}.
\end{equation}

\

Similarly, we multiply \eqref{E:laplace} by $(\mathbf{x},y)\cdot\nabla\phi$ and integrate over $\Omega$
to show that 
\begin{align}\label{E:integral}
0=&\iint_\Omega(\mathbf{x},y)\cdot\nabla\phi\;\Delta\phi~d\mathbf{x}dy \\
=&\int_\Gamma(\mathbf{x}\cdot\nabla_\mathbf{x}\phi+y\phi_y)(\nabla\phi\cdot\boldsymbol{\nu})~dS
-\iint_\Omega\Big(|\nabla\phi|^2+(\mathbf{x},y)\cdot\nabla\Big(\frac12|\nabla\phi|^2\Big)\Big)~d\mathbf{x}dy 
\notag \\
=&-\int_{\mathbb{R}^{n-1}}(\mathbf{x}\cdot[\nabla_\mathbf{x}\phi]_\Gamma
+\eta[\phi_y]_\Gamma)(\mathbf{c}\cdot\nabla\eta)~d\mathbf{x}
-\iint_\Omega|\nabla\phi|^2~d\mathbf{x}dy \notag \\
&-\int_\Gamma\frac12|\nabla\phi|^2(\mathbf{x},y)\cdot\boldsymbol{\nu}~dS
+\iint_\Omega\nabla(\mathbf{x},y)\,\frac12|\nabla\phi|^2~d\mathbf{x}dy \notag \\
=&-\int_{\mathbb{R}^{n-1}}(\mathbf{x}\cdot[\nabla_\mathbf{x}\phi]_\Gamma
+\eta[\phi_y]_\Gamma)(\mathbf{c}\cdot\nabla\eta)~d\mathbf{x} \notag \\
&+\int_{\mathbb{R}^{n-1}}(g\eta-\mathbf{c}\cdot[\nabla\phi]_\Gamma)
(\eta-\mathbf{x}\cdot\nabla\eta)~d\mathbf{x}
+\Big(\frac{n}{2}-1\Big)\iint_\Omega|\nabla\phi|^2~d\mathbf{x}dy \notag \\
=&-\int_{\mathbb{R}^{n-1}}((\mathbf{x}\cdot[\nabla_\mathbf{x}\phi]_\Gamma)(\mathbf{c}\cdot\nabla\eta)
-(\mathbf{c}\cdot[\nabla_\mathbf{x}\phi]_\Gamma)(\mathbf{x}\cdot\nabla\eta))~d\mathbf{x} \notag \\
&-\int_{\mathbb{R}^{n-1}} (\eta[\phi_y]_\Gamma(\mathbf{c}\cdot\nabla\eta)
+\eta\mathbf{c}\cdot[\nabla\phi]_\Gamma)~d\mathbf{x} \notag \\
&+\int_{\mathbb{R}^{n-1}} g\eta(\eta-\mathbf{x}\cdot\nabla\eta)~d\mathbf{x}
+\Big(\frac{n}{2}-1\Big)\iint_\Omega|\nabla\phi|^2~d\mathbf{x}dy=:I_1+I_2+I_3+I_4.\notag
\end{align}
The first equality uses \eqref{E:laplace}, the second equality uses the divergence theorem,
and the third equality uses \eqref{E:kinematic} and the divergence theorem. 
Indeed $[\phi_y]=([\nabla_\mathbf{x}\phi]-\mathbf{c})\cdot\nabla\eta$.
The fourth equality uses \eqref{E:bernoulli}, and the fifth equality uses that $\mathbf{c}\in\mathbb{R}^{n-1}$.
Indeed $\mathbf{c}\cdot[\nabla \phi]_\Gamma=\mathbf{c}\cdot[\nabla_\mathbf{x}\phi]_\Gamma$.

Note from \eqref{E:pohozaev1} that $I_2=-\iint_\Omega|\nabla\phi|^2~d\mathbf{x}dy$. 
Moreover an integration by parts leads to that
\[
I_3=\int_{\mathbb{R}^{n-1}}g\eta^2-g\mathbf{x}\cdot\nabla\Big(\frac12\eta^2\Big)~d\mathbf{x}
=\Big(1+\frac{n-1}{2}\Big)\int_{\mathbb{R}^{n-1}}g\eta^2~d\mathbf{x}.
\]

To proceed, we may assume, without loss of generality, that $\mathbf{c}=(c,0,\dots,0)$. 
We compute that
\begin{align*}
I_1=&-c\int_{\mathbb{R}^{n-1}} \sum_{j=1}^{n-1}
(x_j\phi_{x_j}(\mathbf{x},\eta(\mathbf{x}))\eta_{x_1}(\mathbf{x})
-x_j\phi_{x_1}(\mathbf{x},\eta(\mathbf{x}))\eta_{x_j}(x))~d\mathbf{x} \\
=&-c\int_{\mathbb{R}^{n-1}} \sum_{j=1}^{n-1}
(x_j(\phi_{x_j}(\mathbf{x},\eta(\mathbf{x}))-\phi_y(\mathbf{x},\eta(\mathbf{x}))\eta_{x_j}(\mathbf{x})
+\phi_y(\mathbf{x},\eta(\mathbf{x}))\eta_{x_j}(\mathbf{x}))\eta_{x_1}(\mathbf{x}) \\
&\hspace*{63pt}
-x_j(\phi_{x_1}(\mathbf{x},\eta(\mathbf{x}))-\phi_y(\mathbf{x},\eta(\mathbf{x}))\eta_{x_1}(\mathbf{x})
+\phi_y(\mathbf{x},\eta(\mathbf{x}))\eta_{x_1}(\mathbf{x}))\eta_{x_j}(\mathbf{x}))~d\mathbf{x} \\
=&-c\int_{\mathbb{R}^{n-1}}\sum_{j=1}^{n-1}
(x_j(\phi(\mathbf{x},\eta(\mathbf{x})))_{x_j}\eta_{x_1}(\mathbf{x})
-x_j(\phi(\mathbf{x},\eta(\mathbf{x})))_{x_1}\eta_{x_j}(\mathbf{x}))~d\mathbf{x} \\
=&(n-1)c\int_{\mathbb{R}^{n-1}}\phi(\mathbf{x},\eta(\mathbf{x}))\eta_{x_1}(\mathbf{x})~d\mathbf{x}
+c\int_{\mathbb{R}^{n-1}}\sum_{j=1}^{n-1}
x_j\phi(\mathbf{x},\eta(\mathbf{x}))\eta_{x_1x_j}(\mathbf{x})~d\mathbf{x} \\
&-c\int_{\mathbb{R}^{n-1}}\phi(\mathbf{x},\eta(\mathbf{x}))\eta_{x_1}(\mathbf{x})~d\mathbf{x}
-c\int_{\mathbb{R}^{n-1}}\sum_{j=1}^{n-1}
x_j\phi(\mathbf{x},\eta(\mathbf{x}))\eta_{x_jx_1}(\mathbf{x})~d\mathbf{x} \\
=&(n-2)\int_{\mathbb{R}^{n-1}}\phi(\mathbf{x},\eta(\mathbf{x}))(c\eta_{x_1}(\mathbf{x}))~d\mathbf{x}
=(2-n)\iint_\Omega|\nabla\phi|^2~d\mathbf{x}dy.
\end{align*}
The fourth equality uses integration by parts and the last inequality uses \eqref{E:pohozaev1}.
This proves \eqref{E:pohozaev}. 

\

In case $T>0$, we may rerun the previous argument and calculate that
\begin{align*}
I_3=&\int_{\mathbb{R}^{n-1}}\Big(g\eta-T\nabla\Big(\frac{\nabla\eta}{\sqrt{1+|\nabla\eta|^2}}\Big)\Big)
(\eta-\mathbf{x}\cdot\nabla\eta)~d\mathbf{x} \\
=&\Big(1+\frac{n-1}{2}\Big)\int_{\mathbb{R}^{n-1}}g\eta^2~d\mathbf{x}\\
&+\int_{\mathbb{R}^{n-1}}T\frac{|\nabla\eta|^2}{\sqrt{1+|\nabla\eta|^2}}~d\mathbf{x}
-\int_{\mathbb{R}^{n-1}}T
\frac{|\nabla\eta|^2+\mathbf{x}\cdot\nabla(\frac12|\nabla\eta|^2)}{\sqrt{1+|\nabla\eta|^2}}~d\mathbf{x}\\
=&\Big(1+\frac{n-1}{2}\Big)\int_{\mathbb{R}^{n-1}}g\eta^2~d\mathbf{x}
-\int_{\mathbb{R}^{n-1}} \frac12T\mathbf{x}\cdot\nabla(\sqrt{1+|\nabla\eta|^2}-1)~d\mathbf{x} \\
=&\Big(1+\frac{n-1}{2}\Big)\int_{\mathbb{R}^{n-1}}g\eta^2~d\mathbf{x}
+\frac12\int_{\mathbb{R}^{n-1}}T(\sqrt{1+|\nabla\eta|^2}-1)~d\mathbf{x}.
\end{align*}
This proves \eqref{E:pohozaevT}. 

\begin{remark*}[Extension to the finite depth case]\rm
Let's assume for simplicity that the fluid bottom is rigid and horizontal. One must replace \eqref{E:infinity} by 
\[
\phi_y=0 \qquad \text{at $\{(\mathbf{x},-d): \mathbf{x}\in\mathbb{R}^{n-1}\}$,}
\]
where $d>0$ denotes the asymptotic fluid depth. Abusing notation, let
\[
\Omega=\{(\mathbf{x},y):\mathbf{x}\in\mathbb{R}^{n-1}, -d<y<\eta(\mathbf{x})\}
\quad\text{and}\quad\Gamma=\{(\mathbf{x},\eta(\mathbf{x})):\mathbf{x}\in\mathbb{R}^{n-1}\}.
\]
We rerun the previous argument to show that 
\eqref{E:pohozaev1} holds whereas \eqref{E:integral} becomes
\[
0=\iint_\Omega(\mathbf{x},y)\cdot\nabla\phi\Delta\phi~d\mathbf{x}dy
=I_1+I_2+I_3+I_4-\frac12d\int_{\mathbb{R}^{n-1}}|\nabla\phi(\mathbf{x},-d)|^2~d\mathbf{x},
\]
where $I_1$ through $I_4$ are the same as those in \eqref{E:integral}. 
Therefore \eqref{E:pohozaevT} becomes
\begin{multline}
\frac{n}{2}\iint_\Omega|\nabla\phi|^2~d\mathbf{x}dy
+\frac12d\int_{\mathbb{R}^{n-1}}|\nabla\phi(\mathbf{x},-d)|^2~d\mathbf{x} \\
=\frac{n+1}{2}g\int_{\mathbb{R}^{n-1}}\eta^2~d\mathbf{x}
+\frac12T\int_{\mathbb{R}^{n-1}}(\sqrt{1+|\nabla\eta|^2}-1)~d\mathbf{x}.
\end{multline}

\

In the periodic wave setting, on the other hand, 
$\phi$ is {\em not} periodic and one may not expect a Pohozaev type identity. 
A non-existence proof based upon duality, instead, is found in \cite{Toland-no}, for instance.
\end{remark*}

\subsection*{Acknowledgment}
The author is supported by the National Science Foundation under the grant CAREER DMS-1352597, 
an Alfred P. Sloan research fellowship, an Arnold O. Beckman research award RB14100 
and a Beckman fellowship of the Center for Advanced Study at the University of Illinois at Urbana-Champaign.
She thanks the anonymous referees for their careful reading of the manuscript and numerous helpful suggestions and references.

\bibliographystyle{amsalpha}
\bibliography{solitarybib}

\end{document}